\input amstex
\input pstricks

\def\E{\Bbb E}
\def\F{\Cal F}

\def\Z{\Bbb Z}

\def\cl{\text{cl}}

\def\eps{\varepsilon}

\def\Ext{\text{\rm Ext}}
\def\Ext{\text{Ext}}

\def\FK{\Cal{FK}}

\def\Inter{\text{\rm Int}}

\def\L{\Lambda}

\def\muplnbeta{\mu_{\L(n)}^{+,\beta}}
\def\msub{m_{\text{sub}}}
\def\msup{m_{\text{sup}}}

\def\ph{\widehat{p}}

\def\set#1{\left\{#1\right\}}
\def\sig #1{\sigma(#1)}

\def\tend{\rightarrow}

\def\eps{\varepsilon}

\def\tend{\rightarrow}

\def\Lt{\smash{\widetilde{\Lambda}}}
\def\Lh{\widehat{\Lambda}}

\def\diam{\text{\rm diam}}

\def\und#1{\underline{#1}}

\def\brond{\smash{B \raise 7.9pt\hbox{$\kern -9.4pt{\scriptstyle \circ}$}}}

\def\Mn{M_{\Lambda(n)}}

\def\mathbf#1{{\bold{#1}}}
\def\mathrm#1{{\text{#1}}}
\def\Probn#1{\Phi_{\L(n)}^{w,p}[#1]}

\def\marginal#1{\strut\setbox0=%
      \vtop{\hsize=15mm
            \fiverm
            \textfont0=\scriptfont0 \textfont1=\scriptfont1
            \textfont2=\scriptfont2 \textfont3=\scriptfont3
            \parindent=0pt\baselineskip=7pt
             \raggedright\rightskip=1em plus 2em\hfuzz=.5em\tolerance=9000
             \overfullrule=0pt
             #1}%
      \dimen0=\ht0
      \advance\dimen0 by \dp\strutbox
      \ht0=0pt\dp0=\dimen0
      \vadjust{\kern-\dimen0\moveleft\wd0\box0}%
      \ignorespaces}
\def\showlabelleft#1{\vadjust{\setbox4=\vtop{\overfullrule=0mm\hsize 20mm
\parindent=0pt\fiverm\baselineskip=9pt
\rightskip=4mm\boxit{#1}}
\hbox{\kern-20mm\smash{\raise .5ex\box4}}}}
\documentstyle{amsppt}
\NoBlackBoxes
\nologo

\newcount\spheadno
\global\spheadno=0
\def\nextspheadno{\global\advance\spheadno by1 \the\spheadno \global\headno=0\global\subheadno=0\global\subsubheadno=0\global\cno=0\global\formno=0}

\def\spheadlabel#1{\edef #1{\the\spheadno}}

\newcount\headno
\global\headno=0
\def\nextheadno{\global\advance\headno by1 \the\headno \global\subheadno=0\global\subsubheadno=0}
\def\headnum{\nextheadno.\ }
\def\headlabel#1{\edef #1{\the\headno}}

\newcount\subheadno
\global\subheadno=0
\def\nextsubheadno{\global\advance\subheadno by 1 \the\headno.\the\subheadno\global\subsubheadno=0}
\def\subheadnum{\nextsubheadno. }
\def\subheadlabel#1{\edef #1{\the\headno.\the\subheadno}}

\newcount\subsubheadno
\global\subsubheadno=0
\def\nextsubsubheadno{\global\advance\subsubheadno by 1 \the\headno.\the\subheadno.\the\subsubheadno}
\def\subsubheadnum{\nextsubsubheadno\ }
\def\subsubheadlabel#1{\edef #1{\the\subsubheadno}}

\newcount\formno
\newcount\cno
\def\nextno{\global\advance\cno by1 \the\cno }
\newcount\figno
\def\nextfigno{\global\advance\figno by1 \the\figno }

\def\figlabel#1{\edef #1{\the\figno}}

\global\formno=0
\def\nextformno{\global\advance\formno by1 \the\formno }

\def\blankpage{\vfil\nopagenumbers\eject\line{}\vfil\eject\global\advance\pageno
by -1\footline={\hss\tenrm\folio\hss}}
\def\boxit#1{\vbox{\hrule\hbox{\vrule
\vbox spread 4pt{\vss\hbox spread 4pt{\hss #1\hss}\vss}%
\vrule}\hrule}}

\def\eqnum{\tag{\nextformno}}
\def\eqlabel#1{%
\edef%
#1{\the\formno}%
}
\def\eqref#1{($#1$)}
\def\procnum{\nextno}
\def\proclabel#1{
\edef #1{\the\cno}}
\def\procref#1{$#1$}
\newcount\refno
\global\refno=0
\def\nextrefno{\global\advance\refno by 1 }
\nextrefno\edef\AlexI{\the\refno}
\nextrefno\edef\AlexII{\the\refno}
\nextrefno\edef\AlexIII{\the\refno}
\nextrefno\edef\AlexIV{\the\refno}
\nextrefno\edef\ACC{\the\refno}
\nextrefno\edef\Bar{\the\refno}
\nextrefno\edef\Bodineau{\the\refno}
\nextrefno\edef\BodineauII{\the\refno}
\nextrefno\edef\BricLeboPfister{\the\refno}
\nextrefno\edef\CamiaNewman{\the\refno}
\nextrefno\edef\Cerf{\the\refno}
\nextrefno\edef\CerfSF{\the\refno}
\nextrefno\edef\CerfMessikh{\the\refno}
\nextrefno\edef\CerfPiszI{\the\refno}
\nextrefno\edef\CerfPiszII{\the\refno}
\nextrefno\edef\CCS{\the\refno}
\nextrefno\edef\CourMes{\the\refno}
\nextrefno\edef\DKS{\the\refno}
\nextrefno\edef\ES{\the\refno}
\nextrefno\edef\Ellis{\the\refno}
\nextrefno\edef\ForKast{\the\refno}
\nextrefno\edef\GrimoII{\the\refno}
\nextrefno\edef\Grim{\the\refno}
\nextrefno\edef\GrimII{\the\refno}
\nextrefno\edef\GrimMarst{\the\refno}
\nextrefno\edef\Hoeff{\the\refno}
\nextrefno\edef\IoffeI{\the\refno}
\nextrefno\edef\IoffeII{\the\refno}
\nextrefno\edef\IoffSchon{\the\refno}
\nextrefno\edef\LaaMessRuiz{\the\refno}
\nextrefno\edef\CoyWu{\the\refno}
\nextrefno\edef\Mes{\the\refno}
\nextrefno\edef\Ons{\the\refno}
\nextrefno\edef\Pfis{\the\refno}
\nextrefno\edef\PfisVelen{\the\refno}
\nextrefno\edef\Pisz{\the\refno}
\nextrefno\edef\Smir{\the\refno}
\nextrefno\edef\Slade{\the\refno}
\nextrefno\edef\SmirWer{\the\refno}

\topmatter
\NoRunningHeads

\title
The 2d-Ising model near criticality: a FK-percolation analysis
\endtitle

\author
R. Cerf$\null^1$, R. J. Messikh$\null^2$
\endauthor
\affil
Universit\'e Paris Sud$\null^1$,  EPFL$\null^2$
\endaffil
\address
Universit\'e Paris-Sud, Laboratoire de math\'ematiques, 91405 Orsay, France.
\endaddress
\email
rcerf\@math.u-psud.fr
\endemail
\address
Ecole Polytechnique Federale de Lausanne, cmos, 1015 Lausanne, Switzerland.
\endaddress

\curraddr 64, Rue de rive, 1260 Nyon, Switzerland
\endcurraddr
\email messikh\@gmail.com
\endemail
\date
26 November 2008
\enddate

\keywords
Large deviations, criticality, phase coexistence
\endkeywords
\subjclass
60F10
\endsubjclass
\abstract
We study the 2d-Ising model defined on finite boxes at temperatures that are below but very close from the critical point.
When the temperature approaches the critical point and the size of the box grows fast enough, we establish large deviations estimates on FK-percolation events that concern the phenomenon of phase coexistence.
\endabstract
\endtopmatter
\document
\magnification 1200

\head \headnum{ Introduction}\endhead
The present paper is a study of the influence of criticality on surface order large deviations.
Surface order large deviations occur in supercritical FK-percolation and hence, by the FK-Potts coupling, in the Potts models at sub-critical temperatures.
Originally, the study of such atypical large deviations and their corresponding Wulff construction has started for two dimensional models: the Ising model \cite{\DKS, \IoffeI, \IoffeII, \IoffSchon, \Pfis, \PfisVelen}, independent Bernoulli percolation \cite{\ACC, \AlexIII} and the random cluster model \cite{\AlexIV}.
The just cited papers rely on a direct study of the contours.
This leads to results that go beyond large deviations and give an extensive understanding of phase coexistence in two dimensions and at fixed temperatures.
In higher dimensions, other techniques had to be used to achieve the Wulff construction, \cite{\Bodineau, \Cerf, \CerfPiszI, \CerfPiszII}.
There, the probabilistic estimates rely on  block coarse graining techniques \cite{\Pisz}.
These coarse graining techniques also found applications in other problems not related to the Wulff construction, for example in the study of the random walk on the infinite percolation cluster \cite{\Bar}.
A two-dimensional version of block coarse graining of Pisztora has been given in \cite{\CourMes}, using weak mixing results of Alexander \cite{\AlexII}.

In all the cited works, the percolation parameter (or the temperature) is kept fixed.
The subject of our work is to understand how surface order large deviations and in particular block coarse graining techniques are influenced by criticality.
In other words, our goal is to apply these coarse graining techniques in a joint limit where not only the blocks size increases but also the temperature approaches the critical point from below.
It turns out that the study of block coarse graining in such a joint limit gives rise to several new problems.
Indeed, ideas that are most natural and understood in the fixed temperature case become tricky when we approach criticality.
This gives rise to questions like: how does the empirical density of the infinite cluster converge when we approach the critical point ? or how does the boundary condition influence the configuration inside the box when exponential decay starts to degenerate ? We address and to a certain extend solve these questions in the special case of the 2d-Ising model.

One may wonder why we limit our self to the particular case of the 2d-Ising model.
Indeed, at fixed temperature, block coarse graining techniques are known to be adequate for the study of all FK-percolation models in all dimensions not smaller than two. But even in the fixed temperature case in dimensions higher than three, block coarse graining techniques are known to work up to the critical point only in the percolation model \cite{\GrimMarst} and for the Ising model \cite{\BodineauII}.
Unfortunately very little is known concerning the critical behavior of these models in dimension greater than two. When the dimension is greater than a certain threshold, many of the critical exponents take their so called mean-field values \cite{\Slade}.
Despite these results, to our knowledge, no information is available on the critical behavior of the surface tension, i.e, the exponential price per unit area for the probability of a large interface of co-dimension one.
Therefore we are limited to the two dimensional case, where two potential candidates are possible: site percolation on the triangular lattice, where a lot of progress has been made in the rigorous justification of critical exponents \cite{\Smir, \SmirWer, \CamiaNewman} and the 2d-Ising model where even more accurate information is available thanks to explicit computations, see \cite{\Mes} and the references therein.
Site percolation model would have been an easier model to tackle and the techniques we use could handle this case with straightforward modifications.
But the analysis of the corresponding Wulff construction is still out of reach. The reason for that is related to the open question number $3$ at the end of \cite{\SmirWer}.
Therefore, we chose to treat the 2d-Ising case and proof enough block estimates which  permit the use of the techniques of \cite{\CerfPiszI} to establish the existence of the Wulff shape near criticality \cite{\CerfMessikh} under certain constrains on the simultaneous limit (thermodynamical and going to the critical point).

\subhead\subheadnum{ Statement of the main results}\endsubhead
Our results concern the FK-measures of parameter $q=2$ on finite boxes $\L(n)=(-n/2, n/2]^2\cap\Z^2$, where $n$ is a positive integer.
We denote by $\FK(p, \Lt(n))$ the set of the partially wired FK-measures on boxes $\Lt(n)=(-6n/10, 6n/10]^2\cap\Z^2$ at percolation parameter $p$.
The use of slightly enlarged boxes $\Lt$ is merely technical.
When $p>p_c=\sqrt{2}/(1+\sqrt{2})$, we denote by $\theta(p)$ the density of the infinite cluster.
In what follows we will say that a cluster $C$ of a box $\L$ is crossing, if $C$ intersects all the faces of the boundary of $\L$.
When $p>p_c$, it is known \cite{\CourMes} that up to large deviations of the order of the linear size of the box $\L(n)$, there exists a crossing cluster. It is also known that with overwhelming probability this crossing cluster has a density close to $\theta$ and that the crossing cluster intersect all the sub-boxes of at least logarithmic size.
Our main results essentially state that this qualitative picture still holds when we approach the critical point and let the boxes grow fast enough.
To formulate our results, we define for every box $\L$ the following events:
$$
U(\L)=\left\{\exists\text{ an open crossing cluster }C^{*}\text{ in }\L\right\}.
$$
Moreover, for $M>0$, we define
$$\eqalign{
R(\L,M)&= U(\L)\cap\left\{\text{ every open path }\gamma
\subset\L\text{ with }\diam(\gamma)\geq M\text{ is in }C^{*}\right\}\cr
&\cap\left\{C^{*}\text{ crosses every sub-box of }\L\text{ with diameter }\geq M\right\},
}$$
where $\diam(\gamma)=\max_{x,y\in\gamma}|x-y|$ with $|\cdot|$ denoting the Euclidean norm.

\proclaim{Theorem \procnum}
Let $n>1$ and $a>5$. There
exist two positive constants $\lambda,c=c(a)$ such that if $p>p_c$ and
$n>c(p-p_c)^{-a}$  then
$$
{\forall \Phi\in\FK(\Lt(n),p)}\qquad
\log\Phi[U(\L(n))^c]\leq-\lambda(p-p_c)n.
$$
Moreover,
if $M$ is such that
$$
\frac{\log n}{\kappa(p-p_c)}<M\leq n,\eqnum
$$\eqlabel{\eqsizem}%
with $\kappa>0$ small enough, then
$$
{\forall \Phi\in\FK(\Lt(n),p)}\qquad
\log\Phi[R(\L(n),M)^c]\leq-\lambda(p-p_c)M.
$$
\endproclaim\proclabel{\thmUR}%
\noindent Note that the speed of the large deviations slows down by a  factor $(p_c-p)$ when $p\downarrow p_c$.
This is directly related to the critical exponent $\nu=1$ of the inverse correlation length of the 2d-Ising model.
The exponent $a>5$ restrict our result to be valid only for boxes of width  much larger than the inverse correlation length.

\noindent Next, we consider deviations for empirical densities of the infinite cluster when $p\downarrow p_c$.
For $n>0$, we consider the number of boundary connected sites
$$
M_{\L(n)}=|\{x\in\L(n):\,x\leftrightarrow\partial\L(n)\},
$$
where we have used the notation $|E|$ to denote the cardinality of a set $E\subset\Z^2$ and where $\partial\L$ denotes the site boundary of $\L$.
It is known that for all $p>p_c$,
$$
\lim_{n\rightarrow\infty}{1\over |\L(n)|}\Phi^{w,p}_{\L(n)}[M_{\L(n)}]=\theta(p).\eqnum
$$\eqlabel{\convtheta}%
On the other hand, from the solution of Onsager \cite{\Ons} we know that $\theta(p)\sim(p-p_c)^{1/8}$ when $p\downarrow p_c$.
This degeneracy requires us to control the speed at which the convergence \eqref{\convtheta} occurs.
To this end, for each  $\delta>0$, we define
$$
\msup(\delta,p)=\inf\left\{m\geq 1:\,\forall n\geq m\quad \Phi_{\L(n)}^{w,p}[M_{\L(n)}]\leq|\L(n)|(1+\delta/2)\theta\right\},
$$
which represents the minimal size of the box required to approximate the density of the infinite cluster within an error of $\delta\theta/2$.
The subadditivity of the map $\L\mapsto M_{\L}$, makes it handy to consider large deviations from above.
To do so, we define the event
$$
W(\L,\delta)=\left\{M_{\L}\leq(1+\delta)\theta|\L|\right\}.
$$
and obtain
\vfill\eject
\proclaim{Theorem \procnum}
Let $p>p_c$ and $\delta>0$.
If $n>8\msup(\delta,p)/\delta$ then
$$
\log\Phi^{w,p}_{\L(n)}[W(\L(n),\delta)^c]\leq-\left(\frac{\delta\theta n}{4\msup(\delta,p)}\right)^2.\eqnum
$$
In particular, for every $a>5/4$, there exists a positive
constant $c=c(a,\delta)$ such that whenever $n\uparrow\infty$ and
$p\downarrow p_c$ in such a way that $n>c(p-p_c)^{-a}$ then
$$
\lim_{n,p}{1\over (p-p_c)^{2a+1/4}n^2}\log\Phi^{w,p}_{\L(n)}[W(\L(n),\delta)^c]\,<0\,.
$$
\eqlabel{\eqIW}%
\endproclaim\proclabel{\thmW}%
\noindent It is natural to take the density of the crossing cluster as an empirical density of the infinite cluster.
Next we consider the deviations from below of this quantity.
For any $\delta >0$, we define the event
$$
V(\L,\delta)=U(\L)\cap\left\{|C^{*}|\geq(1-\delta)\theta|\L|\right\}.
$$
When $p>p_c$ is kept fixed, an upper bound of the correct exponential speed can be obtained using coarse graining techniques of Pisztora.
We proof that similar ideas can be used to obtain a priori estimates in the joint limit.
\proclaim{Theorem \procnum}
Let $a>5$ and $\alpha\in ]0,(1+{1\over 8a})^{-1}[$.
There exists a positive constant
$c=c(a,\alpha)$ such that, if $n\uparrow\infty$ and $p\downarrow p_c$ in
such a way that
$n^{\alpha}(p-p_c)^{a}>c$ then
$$
\sup_{\Phi\in\FK(\Lt(n),p)}\Phi[V(\L(n),\delta)^c]\leq\exp(-\lambda\delta(p-p_c)n^\alpha)+\exp(-{\delta^2\over 4}(p-p_c)^{1/4}n^{2-2\alpha}),\eqnum
$$\eqlabel{\bV}%
where $\lambda$ is a positive constant. In particular
$$
\lim_{n,p}\inf_{\Phi\in\FK(\Lt(n),p)}\Phi[V(\L(n),\delta)]=1.
$$
\endproclaim\proclabel{\thmV}%
\noindent When $p>p_c$ is kept fixed, the right hand side of \eqref{\bV} can be replaced by an expression of the form $\exp(-cn)$ where $c$ is a positive constant.
The appearance of two terms in the joint limit $n\rightarrow\infty$ and $p\downarrow p_c$ comes from the fact that the size of the blocks in the coarse graining cannot be taken constant anymore, they have to diverge like $n^\alpha$.
Note that the two terms on the right hand side of  \eqref{\bV} are competing, indeed when $\alpha$ increases the first term decreases and the second one increases.
\subhead{\subheadnum Organisation of the paper}\endsubhead
In section 2, we start by introducing the basic definitions and notations used in the rest of the paper. In this section, we also provide preliminary results on the critical behavior of the 2d-Ising model.
Then, in section 3, we establish weak mixing results in a situation where $p\rightarrow p_c$. These results will enable us to control adequately the influence of the boundary conditions. Finally, the proofs of the main theorem are given in section 4.  In the appendix, we prove a technical result concerning the speed of convergence of the empirical magnetization near criticality.
\head\headnum Preliminaries\endhead
\subhead\subheadnum{The FK-representation}\endsubhead\subheadlabel{\subhdFKrep}%
There exists a useful and well known coupling between the Ising model
at inverse temperature $\beta$
and the random cluster model with parameter $q=2$ and
$p=1-\exp(-2\beta)$, see \cite{\ES, \ForKast}.
The coupling is a probability
\def\Coupling{\Bbb P_n^{+}}measure $\Coupling$ on the edge-spin configuration space
$\{0,1\}^{\E(\L(n))}\times\{-1,+1\}^{\L(n)}$.

To construct $\Coupling$ we first consider Bernoulli percolation of parameter $p$
on the edge space $\{0,1\}^{\E(\L(n))}$, then we choose the spins of the sites in $\L(n)$ independently with the
uniform distribution on $\{-1,+1\}$ and finally  we condition the edge-spin
configuration on the event that there is no open edge in $\L(n)$ between
two sites with different spin values. The construction can be summed up with a formula, we have
$$\eqalign{
\forall (\sigma,\omega)&\in\{0,1\}^{\E(\L(n))}\times\{-1,+1\}^{\L(n)}\cr
&\Coupling(\sigma,\omega)={1\over Z}\prod_{e\in\E(\L(n))}p^{\omega(e)}(1-p)^{1-\omega(e)}1_{(\sigma(x)-\sigma(y))\omega(e)=0},
}$$
\noindent where $Z$ is the appropriate normalization factor. It can be verified
that the marginal of $\Coupling$ on the spin configurations is the
Ising model at inverse temperature  $\beta$ given by the formula
$p=1-\exp(-2\beta)$ and the marginal on the edge configurations is the
random cluster measure with parameters $p$, $q=2$ and subject to
wired boundary conditions, i.e., the
probability measure on $\Omega_{\L(n)}=\{0,1\}^{\E(\L(n))}$ defined by
$$
\forall\omega\in
\Omega_{\L(n)}\qquad
\Phi_{\L(n)}^{p,w}[\omega]={1\over
Z}q^{\cl^w(\omega)}\prod_{e\in\E(\L(n))}p^{\omega(e)}(1-p)^{1-\omega(e)},\eqnum
$$
\eqlabel{\eqFKmeas}where $\cl^w(\omega)$ is the number of connected components with the
  convention that two clusters that touch the boundary
$\partial\L(n)$ are identified.
  This coupling says that one may obtain an Ising
configuration by first drawing a FK-percolation configuration with the
measure
$\smash{\Phi_{\L(n)}^{w,p}}$, then coloring all the sites in the clusters that
touch the boundary $\partial\L(n)$ in $+1$ and finally coloring the remaining
clusters independently  in $+1$ and $-1$
with probability $1/2$ each. Also, the coupling
permits to obtain a $\smash{\Phi_{\L(n)}^{w,p}}$ percolation configuration by
first drawing a spin configuration with $\muplnbeta$, then  declaring
that all the edges between two sites with different spins are closed,
while
the other edges are independently declared open with probability
$p$ and closed with probability $1-p$.
\par\noindent Let $\L\subset\Z^2$ and $0\leq p\leq 1$. In addition to the wired boundary conditions we
  will also work with {\it partially wired boundary conditions}. In order to
  define them, we consider a partition $\pi$ of $\partial\L=\{x\in\L:
  \exists y\in\Z^2\setminus\L,\,|x-y|_1=1\}$. Let us say that $\pi$ consists of
  $\{B_1,\cdots,B_k\}$, where the $B_i$ are non-empty disjoint subsets of
  $\partial\L$ and such that $\cup_i B_i=\partial\L$. For every
  configuration $\omega\in\Omega_\L$, we define $\cl^{\pi}(\omega)$ as
  the number of open connected clusters in $\L$ computed by
  identifying two clusters that are connected to the same set $B_i$. The $\pi$-wired
  FK-measure $\Phi^{p,\pi}_{\L}$ is defined by substituting
  $\cl^w(\omega)$ for  $\cl^{\pi}(\omega)$ in \eqref{\eqFKmeas}. We will
  denote the set of all partially wired FK-measures in $\L$ by
  $\FK(p,\L)$. Note that $\Phi^{p,w}_\L$ corresponds to
  $\pi=\{\partial\L\}$. We define the FK-measure with free boundary
  conditions $\Phi^{p,f}_\L$ as the partially wired measure
  corresponding to $\pi=\emptyset$.
\par\noindent Let $U\subseteq V\subseteq\Z^2$. For every configuration
  $\omega\in\{0,1\}^{\E(\Z^2)}$, we denote by $\omega_V$ the
  restriction of $\omega$ to $\Omega_V=\{0,1\}^{\E(V)}$. More
  generally we will denote by $\omega_V^U$ the restriction of
  $\omega$ to $\Omega_V^U=\{0,1\}^{\E(V)\setminus\E(U)}$. If $V=\Z^2$
  or $U=\emptyset$ then we drop them from the notation. We will
  denote by $\F_V^U$ the  $\sigma$-algebra generated by the finite
  dimensional cylinders of $\Omega_V^U$.

\par\noindent Note that every configuration $\eta\in\Omega_V$ induces
  a partially wired boundary condition $\pi(\eta)$ on the set $U$. The partition
  $\pi(\eta)$ is obtained by identifying the sites of $\partial U$
  that are connected through an open path of  $\eta^U$. We will denote
  by $\Phi^{p,\pi(\eta)}_U$ the corresponding FK measure.
\subhead\subheadnum{ Planar Duality}\endsubhead The duality of the FK-measures in dimension two
  is well known. In this paper we will use the notation of
  \cite{\CourMes} that we summarize next. Let $0\leq p\leq 1$ and $\L$
  be a box of $\Z^2$.
To construct the dual model we associate to a box ${\L}$ the set
 $\widehat{{\L}}\subset\Z^2+(1/2,1/2)$, which is defined as the
smallest box of $\Z^2+(1/2,1/2)$ containing
${\L}$.
To each edge $e\in\E({\L})$ we associate the edge
 $\widehat{e}\in\E(\widehat{{\L}})$ that crosses the edge
 $e$. Note that $\{e'\in\E(\widehat{{\L}}):\exists
 e\in\E({\L}),\widehat{e}=e'\}=\E(\widehat{{\L}})\setminus\E(\partial\widehat{{\L}})$.

\noindent This allows us to build a bijective application from
$\Omega_{{\L}}$ to $\Omega_{\widehat{{\L}}}^{\partial\widehat{{\L}}}$ that maps each original configuration
 $\omega\in\Omega_{{\L}}$ into its dual configuration $\widehat{\omega}
 \in\Omega_{\widehat{{\L}}}^{\partial\widehat{{\L}}}$ such
 that
$$
\forall e\in\E({\L}):
 \widehat{\omega}(\widehat{e})=1-\omega(e).
$$
The duality property states that for any $0<p<1$ and any  $\Cal
F_{\L}$-measurable event $A$  we have
$$
\Phi^{f,p}_{{\L}}[A]=\Phi^{w,\widehat{p}}_{\widehat{{\L}}}[\widehat{A}],
$$
where $\widehat{A}=\{\eta\in\Omega_{\widehat{{\L}}}: \exists
\omega\in A, \widehat{\omega}=\eta^{\partial{\widehat{{\L}}}}\}\subset
\Omega_{\widehat{{\L}}}^{\partial\widehat{{\L}}}$ is the dual event of $A$
and where $\widehat p=2(1-p)/(2-p)$.
It is useful to remark that when we translate an $\Cal F_{\L}$-measurable event $A$ into
it's dual  $\widehat{A}$, we obtain an event which is in $\Cal
F_{\widehat{{\L}}}^{\partial\widehat{{\L}}}$. and that $\Phi^{w,\widehat{p},q}_{\widehat{{\L}}}[\widehat{A}]$ does note depend on  the states of the edges in $\E(\partial\widehat{{\L}})$.
Note also that under the measure $\Phi^{w,p}_\L$ the law of $\omega_{\partial\L}$ is an independent percolation of parameter $p$ and $\omega_{\partial\L}$ is also independent from $\omega^{\partial\L}$.

\noindent We end this section by setting the following convention concerning the use of the word {\it dual} in the rest of the paper: we always consider that the original  model is the super-critical one, i.e.,  $p>p_c$, which is defined on the edges of $\Z^2$. The {\it dual } model is always the dual of the super-critical model.
That is, it is a sub-critical model defined on the edges of $\Z^2+(1/2,1/2)$ and at percolation parameter $\widehat p=2(1-p)/(2-p)\leq p_c$.
A dual path, circuit or site will always denote a path, circuit or site in $\Z^2+(1/2,1/2)$.
The term {\it open dual} will always designate edges $\widehat e$ of $\Z^2+(1/2,1/2)$ that are open with respect to the dual configuration,
i.e., $\widehat\omega(\widehat e)=1$.
The law of the dual edges $\widehat e$  will always be the dual measure $\Phi^{\widehat p}$ which is sub-critical, i.e., $\widehat p < p_c$.

\subhead\subheadnum{Preliminary results on criticality in the 2d-Ising model}\endsubhead
In this section we review some known results about the nature of the phase transition of the 2d-Ising model.
These properties are important for our analysis and their proofs uses the specificities of the 2d-Ising model: explicit computations and correlation inequalities.
Even though similar results are believed to hold for all the two dimensional FK-measures with parameter $1\leq q\leq 4$, the FK-percolation with parameter $q=2$ is the only model where such results can be established via-explicit computations. That is why our results are restricted to the 2d-Ising model.
Let us also mention that if the analogues of the results stated in the section where available for other two dimensional FK-measures then the techniques used in this paper can be generalized to treat such cases. The extension to higher dimensional models is potentially also possible along the ideas of \cite{\Pisz} but, to our knowledge,  information about the critical behavior of the surface tension near criticality is nowadays unavailable even in the form of conjectures.
\subsubhead\subsubheadnum{ The critical point}\endsubsubhead It is known that the critical point of the Ising
model on $\Z^2$ is given by the fixed point of a duality relation (see \cite{\Grim}). For the
random cluster model with $q=2$, the dual point $\ph$ is related to
$p$ through the relation
$$
{p\over 1-p}\ {\ph\over 1-\ph}=2,
\text{ and the fixed point is }p_c={\sqrt 2\over1+\sqrt{2}}.\eqnum
$$
\eqlabel{\eqrelationduality}%

\noindent For the general $q$-Potts model, the identification of the critical point and the self-dual point, i.e.,
$p_c=\sqrt{q}/(1+\sqrt{q}) $, is still an open problem for the values
$2<q<25$. When $q>25.72$, this identity has been established and in
this situation the Potts model exhibits a first order phase transition
\cite{\GrimoII,\LaaMessRuiz}. Thus the 2d-Ising model is the only
two  dimensional Potts model exhibiting a second order phase transition for
which the critical point has been rigorously identified to be the self-dual point.

\subsubhead\subsubheadnum{The surface tension}\endsubsubhead
In the two dimensional supercritical FK-percolation model, large interfaces are best studied via duality.
Indeed, a large interface implies a long connection in the sub-critical  dual model.
This is why the surface tension at $p>p_c$ is given by the exponential decay of connectivities in the sub-critical dual model:
$$
\forall x\in\Z^2\quad \tau_p(x)=-\lim_{n\tend\infty}{1\over
n}\log\Phi^{\widehat p}_\infty[0\leftrightarrow nx],
$$
where $\Phi^{\widehat p}_\infty$ denotes the unique infinite FK-measure for $\widehat p<p_c$ \cite{\GrimII}.
In this paper, we are interested in the situation where the spatial
scale $n$ goes to infinity and simultaneously $p$ goes to
$p_c$.
Using sub-additivity and the formula for $\tau_p$, it is possible to show that
\proclaim{Proposition \procnum} When
$n\uparrow\infty$ and $p\downarrow p_c$ we have uniformly in $x\in Z^2$ that
$$
{1\over (p-p_c)n|x|}\log\Phi^{\widehat p}_\infty[0\leftrightarrow nx]\leq-\tau_c,\eqnum
$$
\eqlabel{\eqdecexplim}where $\tau_c$ is a positive constant.
\endproclaim
\proclabel{\procLimiteJointe}
\noindent The proof of the last proposition and even stronger results is the subject of \cite{\Mes}.

\subsubhead\subsubheadnum{The magnetization}\endsubsubhead The magnetization of the Ising model corresponds to the density $\theta(p)$ of the infinite cluster in the FK-representation.
When $q=2$, it is known that $\theta(p)$ approaches zero when $p\downarrow p_c$.
Thanks to the Onsager's exact solution, it is also known at which speed this occurs:
$$
\theta(p)\sim(p-p_c)^{1/8}\quad\text{when}\, p\downarrow p_c. \eqnum
$$\eqlabel{\eqOnsager}%
To apply our techniques, we will also need to know at which speed the empirical magnetization converges to $\theta(p)$ when approaching $p_c$.
More precisely, we need to control
$$
{1\over n^2}\Phi^{p,w}\left [|\{x\in\L(n):\,x\leftrightarrow\partial\L(n)\}|\right]-\theta(p)\eqnum
$$\eqlabel{\speedemp}%
in the joint limit $n\rightarrow\infty$ and $p\rightarrow p_c$.
In turns out that the control of \eqref{\speedemp} is delicate.
Indeed, we where unable to control the speed of convergence of \eqref{\speedemp} in the joint limit using only  Proposition \procref{\procLimiteJointe}, \eqref{\eqOnsager} and robust FK-percolation techniques.
We found a solution to this problem using further specificities of the 2d-Ising model, namely correlation inequalities.
Using the ideas of \cite{\BricLeboPfister}, we get the following result
\proclaim{Proposition \procnum}Let $\xi>0$ and $a>\xi+1$. There
exist two positive constants $c=c(\xi,a)$ and $\rho$ such that
$$
\forall p\neq p_c,\,n>c|p-p_c|^{-a}\qquad{1\over n^2}\Phi^{p,w}\left [\sum_{x\in\L(n)}1_{x\leftrightarrow\partial\L(n)}\right]-\theta(p)\leq \,\rho|p-p_c|^\xi.
$$
\endproclaim\proclabel{\procregpfister}%
\noindent We defer the proof of the last proposition to the end of the paper in Appendix A.

\head\headnum{ Weak mixing near criticality}\endhead
In this part we establish weak mixing properties in
the situation where $p\downarrow p_c$. These results are crucial in order to bound the influence
of the boundary conditions. As it appears from \cite{\CourMes}, in order to
implement a useful coarse graining in dimension two, it is necessary to
have a control of the boundary conditions. When $p$ is
fixed, this control can be obtained by using the weak mixing
properties proved in \cite{\AlexI, \AlexII}. To handle the situation
where $p\downarrow p_c$, we give an alternative way to establish weak mixing and generalize the results of
\cite{\AlexI, \AlexII} to a situation where the exponential decay of connectivities
becomes degenerate.

\subhead\subheadnum{Control of the number of boundary connected sites}\endsubhead Let
$p<p_c, n\geq 1$. In this paragraph, we are
interested in the control of the number of boundary connected sites
$$
\Mn=|\set{x\in \Lambda(n): x\leftrightarrow\partial\Lambda(n)}|.\eqnum
$$
\eqlabel{\eqdefbdvertices}The coming results depend on the speed of convergence of
the mean of $M_n$ near the critical point. We characterize this speed
by introducing the following quantity:
$$
\forall p<p_c,\,\delta>0\quad\msub(\delta,p)=\inf\left\{m\geq 1:\, \forall n>m
\quad {1\over|\L(n)|}\Phi^{w,p}_{\L(n)}[M_{\L(n)}]\leq\delta\right\}.\eqnum
$$\eqlabel{\eqdefmsub}

The main tool used in this section is subadditivity which permit us to reduce the problem to a family of bounded i.i.d random variables.
Which are then well under control thanks to the following concentration bound:
\proclaim{Lemma \procnum}(Theorem 1 of \cite{\Hoeff}) If $(X_i)_{1\leq
i\leq n}$ are independent random variables
with values in $[0,1]$ and with mean $m$, then
$$
\forall t\in ]0,1-m[\qquad
P\Big[\sum_{i=1}^{n}\left(X_i-m\right)\geq n\ t
\Big]\leq\exp(-nt^2).
$$
\endproclaim\proclabel{\BornHoeff}

\proclaim{Lemma \procnum}
Let $\delta>0,\ p\leq p_c$. If
$
n\geq16\msub(\delta/2,p)/\delta,
$
then
$$
\log\Phi^{w,p}_{\Lambda(n)}\left[\frac{\Mn}{|\L(n)|}\geq\delta\right]\leq-\left(\frac{\delta
n}{6\msub(\delta/2,p)}\right)^2.
$$
\endproclaim\proclabel{\procbern}
\demo{Proof}
First we partition  $\Lambda(n)$ into
translates of the square $\Lambda(m)$ where
$$
m=\msub(\delta/2,p).\eqnum
$$
\eqlabel{\eqSizemVarphi}Next, we take
$$
n>16m/\delta,\eqnum
$$
\eqlabel{\eqSizemN}
\noindent and consider the set
$$
\L'(n)=\bigcup_{\und x\in\Z^2:B(\und x)\subset\L(n)}B(\und x),
$$
\noindent where $B(\und x)=m\und x+\L(m)$.
Note that $|\L(n)\setminus\L'(n)|\leq 4mn$. The number
of
partitioning blocks  satisfies
$$
\frac{n^2}{2m^2}\leq|\und\L'(n)|\leq \frac{n^2}{m^2}.\eqnum
$$
\eqlabel{\eqBorninftaille}
\noindent Since $M_\Lambda$ is subadditive, by
\eqref{\eqBorninftaille} and \eqref{\eqSizemN}, we obtain
$$\eqalign{
\frac{\Mn}{|\L(n)|}&\leq\frac{1}{n^2}\sum_{\und x\in\und \L'(n)}|\{v\in
B(\und x): v\leftrightarrow\partial\L(n)\}|+\frac{4m}{n}
\cr
&\leq\frac{1}{|\und\L'(n)|}\sum_{\und x\in\und \L'(n)}\frac{M_{B(\und x)}}{|B(\und x)|}+\frac{\delta}{4}.
}
$$By the FKG inequality, we get
$$
\Phi^{w,p}_{\Lambda(n)}\left[\frac{\Mn}{|\L(n)|}\geq\delta\right]\leq\Phi^{w,p}_{\Lambda(n)}\left[\frac{1}{|\und\L'(n)|}\left.\sum_{\und
x\in\und \L'(n)}
\frac{M_{B(\und x)}}{|B(\und x)|}\geq\frac{3\delta}{4}\right\vert
E\right]\eqnum
$$
\eqlabel{\eqrecentring}where $E$ is the increasing event
$\left\{\forall\und x\in\und\L'(n),\ \text{all the edges of } \partial B(\und x)
\text{ are open}\right\}$.
The random variables $M_{B(\und x)}/|B(\und x)|$,
$\und x\in \und\L'(n)$,
take their values in $[0,1]$ and they are
independent under $\Phi_{\L(n)}^{w,p}[\cdot\ \vert E]$.
By \eqref{\eqSizemVarphi}, their mean satisfies
$$
\forall\ \und x\in \und\L'(n)\qquad\Phi_{\L(n)}^{w,p}\left.\left[\frac{M_{B(\und x)}}{|B(\und
x)|}\right\vert E\right]=\Phi_{B(\und x)}^{w,p}\left[\frac{M_{B(\und
x)}}{|B(\und x)|}\right]\leq\frac{\delta}{2}.\eqnum
$$
\eqlabel{\eqBornMoy}Finally, by lemma \procref{\BornHoeff} and by the
inequalities \eqref{\eqBorninftaille}, \eqref{\eqrecentring} and \eqref{\eqBornMoy}  we get
$$
\eqalign{
\Phi^{w,p}_{\Lambda(n)}\left[\frac{\Mn}{|\L(n)|}\geq\delta\right]&\leq\exp\left(-\frac{\delta^2n^2}{32m^2}\right).
}
$$
\qed
\enddemo
\subhead\subheadnum{ Control of the boundary conditions }\endsubhead
In this section, we determine a regime where we can still control the
influence of the boundary conditions when $p\tend p_c$. The regime will be characterized by the speed by which the quantity $\msub$
defined in \eqref{\eqdefmsub} diverges near the critical point.
We thus need to give an upper bound for the speed of this divergence.

\proclaim{Lemma \procnum} Let $\kappa>0, \xi>0$. For every $a>\xi+1$ there exists
a positive constant $c=c(a,\kappa)$  such that
$$
\forall p<p_c\qquad
\msub(\kappa(p_c-p)^{\xi},p)\leq c(p_c-p)^{-a}.
$$
\endproclaim\proclabel{\procSpeedmsub}
\demo{Proof}Let $a>1$ and $\xi\in(0,a-1)$. From
proposition \procref{\procregpfister} we know that for every
$\eta\in(\xi,\xi+1)$ there exist two positive
constants $\rho$ and $c_1$ such that
$$
\forall p<p_c\quad
\forall\,n>c_1(p_c-p)^{-a}\qquad{1\over|\L(n)|}\Phi^{w,p}_{\L(n)}[M_{\L(n)}]\leq \,\rho(p_c-p)^\eta.
$$Furthermore, since $\eta>\xi$, there exists a positive constant $\eps=\eps(\rho,\xi,\kappa,\eta)$
such that
$$
\forall\,p\in(p_c-\eps,p_c)\qquad \rho(p_c-p)^\eta\leq\kappa(p_c-p)^{\xi}.
$$Note also that if $p\leq p_c-\eps$ then $\kappa(p-p_c)^\xi\geq\kappa\eps^\xi$ and there
exists $n_0(\eps^\xi)$  such that $n>n_0$ implies
$$
\forall\,p<p_c-\eps\qquad\msub(\kappa(p-p_c)^\xi,p)\leq n_0.
$$Hence the result follows by choosing
$c=\max(c_1,\eps^an_0)$. \qed
\enddemo
\proclaim{Proposition \procnum} Let $p<p_c$ and $a>5$. There exist two positive
constants $c=c(a)$ and $\lambda$  such that if $n>c(p_c-p)^{-a}$ then
$$
\log\Phi_{\L(n)}^{w,p}[0\leftrightarrow\partial\L(n)]\leq -\lambda(p_c-p)n
\,.
$$
\endproclaim\proclabel{\procContConn}
\demo{Proof}
Let $A=\{0\leftrightarrow\partial\L(n/2)\}$. In order to
control the influence of the boundary conditions imposed on $\L(n)$
we first write
$$\eqalign{
\Phi_{\Lambda(n)}^{w,p}[A]\leq&\Phi_{\Lambda(n)}^{w,p}[A\cap\{\Mn\leq|\L(n)|\delta\}]\cr
&+\Phi_{\Lambda(n)}^{w,p}[\Mn>|\L(n)|\delta\ ],
}\eqnum$$
\eqlabel{\eqfirstdecomp}where $\Mn$ is defined in \eqref{\eqdefbdvertices}.
On the event $A'=A\cap\{\Mn\leq|\L(n)|\delta\}$ of the first term we can bound the influence of the
boundary conditions in an adequate way
by using a judicious trick due to David Barbato \cite{\Bar},
while the second term will be
made negligible thanks to lemma \procref{\procbern}.

\noindent{\it Barbato's trick:} This trick has initially been
introduced in order to simplify the proof of the so called interface
lemma in the case of dimensions higher or equal to three. Here we will
use this trick in a different context. From the definition of the FK-measures it is clear that the influence
of the boundary conditions comes from the connected components that
connect $\partial\L(n/2)$ to $\partial\L(n)$. Thus if one can
cut all these connections without altering too much the probability of
the event $A$ then one gets a control over the influence of the
boundary conditions. To do this we first  define
$M'_{\L(n)}$ as
$$M'_{\L(n)}\,=\,
\Big|\{x\in\L(n):
x\leftrightarrow\partial\L(n)\text{ in }\L(n)\setminus\L(2|x|_\infty)\}
\Big|\,.$$This is the same quantity as $\Mn$ with the difference that we count
only the sites $x$ that are connected to the boundary with a {\it direct} path
that does not use the edges in $\E(\L(2|x|_\infty))$.
Now suppose that
$A'=A\cap\{M_{\L(n)}\leq|\L(n)|\delta\}$ occurs.
Since
$M'_{\L(n)}\leq\Mn$ we also have
$M'_{\L(n)}\leq\delta|\L(n)|$. Next, for $0<h<1/4$, we define the set
$$
\frak b(h)=\partial[-n(1-h)/2,n(1-h)/2]^2.
$$Note that for $0<h<1/4$, we always have
$$
\frak b(h)\cap\L(n/2)=\emptyset.
$$Next, we concentrate on the finite set of values $0<h_1<\cdots<h_K$
that satisfy
$$
\frak b(h_k)\cap\L(n)\neq\emptyset.
$$We notice that the number $K$ of such values $h_k$ satisfies
$$
{n\over 8}-1<K<{n\over 8}+1.
$$
\noindent Until here, the construction does not depend on the configuration. Next, we scan
the configuration in $\L(n)$ from outside inwards and define for each
$h_k$ the set of {\sl bad} sites intersected by $\frak b(h_k)$:
$$
V(h_k)\,=\,M'_{\L(n)}\cap\frak b(h_k).
$$On $A'$ we have that $\sum_{k=1}^{K}|V(h_k)|\leq M'_{\Lambda(n)}\leq\delta|\Lambda(n)|$
whence, for $n$ large enough,
$$
\min_{k}|V(h_k)|\leq{\delta|\Lambda(n)|\over K}\leq{\delta |\Lambda(n)|\over {n\over 8}-1}\,\leq\,16\delta n\,.
$$Thus there exists at least one $k\in\{1,\dots,K\}$ such that
$$
|V(h_k)|\leq
16\delta n
.\eqnum
$$
\eqlabel{\eqbadvertices}\def\Carre#1{\L((1-#1)n)}We define $h^*$ as the first (smallest) value $h_k$ that
satisfies \eqref{\eqbadvertices}. Notice that $h^{*}$ is a sort of stopping time, in the sense that
$$
\forall 0<h<1/4\qquad \{h^{*}=h\}\in\F_{\L(n)\setminus\Carre{h}}.\eqnum
$$
\eqlabel{\eqStopTimeh}Then we define the set of {\sl bad} edges as the set of edges that
have one extremity in $\Carre{h^*}$ and the other in $V(h^*)$:
$$
I_n\,=\,
\big\{\,e=\{v,u\}\in\E^2:\ v\in\Carre{h^{*}},\ u\in V(h^{*})\,
\big\}.
$$Even though
$$
I_n\cap\E(\L(n)\setminus\Carre{h^*})=\emptyset,\eqnum
$$
\eqlabel{\eqInsideI}we obtain from \eqref{\eqStopTimeh} and from the
definition of $V(h^*)$ that
$$
\forall I\subseteq\E(\L(n))\quad\{I_n=I\}\in\F_{\L(n)\setminus\Carre{h^{*}}}.\eqnum
$$
\eqlabel{\eqOutsidemeasI}It is also important to notice that
$$
I_n\cap\E(\L(n/2))=\emptyset.\eqnum
$$
\eqlabel{\eqOutsideI}Now, for each site $v\in V(h^*)$ there is at
most one edge $e$ in $I_n$ with
extremity $v$ thus we get from \eqref{\eqbadvertices}  that
$$
|I_n|\leq 16\delta n.\eqnum
$$
\eqlabel{\eqbadedges}Let $\Psi:A'\rightarrow\Omega$ be
the map defined by:
$$
\forall\omega\in A'\quad\forall e\in\L(n)\quad
\Psi(\omega)(e)=\left\{\eqalign{0&\qquad\text{if }e\in
I_n(\omega)\cr\omega(&e)\quad\text{otherwise}\cr}\right.
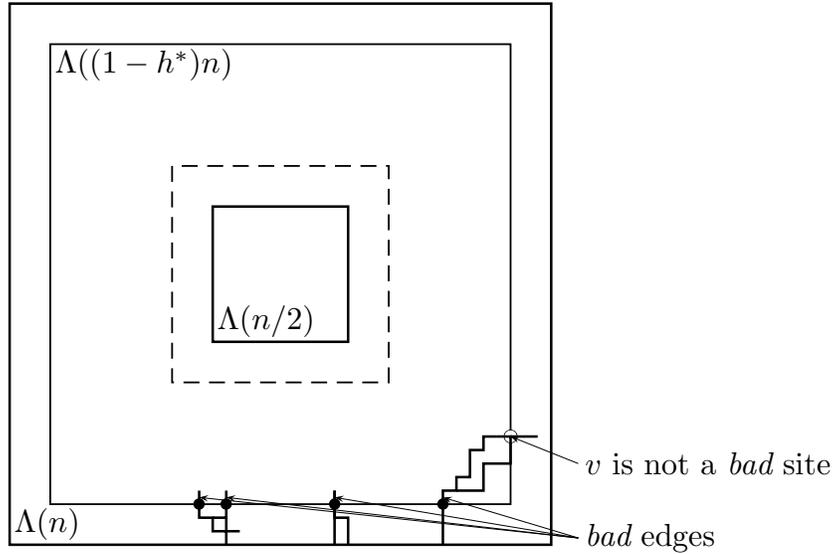
$$

\vbox{\centerline{
\psset{unit=0.15cm}
\pspicture(-20,-25)(20,20)
\uput{0.5}[45](-20,-20){$\L(n)$}
\pspolygon(-20,-20)(20,-20)(20,20)(-20,20)
\uput{0.5}[45](-5,-5){$\L(n/2)$}
\pspolygon(-5,-5)(5,-5)(5,5)(-5,5)
\pspolygon[linestyle=dashed,linewidth=0.02cm](-8,-8)(8,-8)(8,8)(-8,8)
\uput{0.5}[-45](-17,17){$\L((1-h^*)n)$}
\pspolygon[linewidth=0.02cm](-17,-17)(17,-17)(17,17)(-17,17)
\psdots[dotscale=1](-6,-17)(-4,-17)(4,-17)(12,-17)
\psline[linewidth=0.01cm]{->}(22,-19.5)(-6,-16.5)
\psline[linewidth=0.01cm]{->}(22,-19.5)(-4,-16.5)
\psline[linewidth=0.01cm]{->}(22,-19.5)(4,-16.5)
\psline[linewidth=0.01cm]{->}(22,-19.5)(12,-16.5)
\uput{0.5}[0](22,-19.5){{\it bad} edges}
\psline(-6,-16)(-6,-18)(-5,-18)(-5,-19)(-4,-19)(-4,-20)
\psline(-5,-18)(-4,-18)(-4,-16)
\psline(-4,-18)(-4,-19)(-3,-19)
\psline(4,-16)(4,-18)(4,-19)(4,-20)
\psline(4,-16)(4,-18)(5,-18)(5,-19)(5,-20)
\psdots[dotscale=1.2,dotstyle=o](17,-12)
\psline(12,-20)(12,-19)(12,-18)(12,-17)
\psline(12,-17)(12,-16)(13,-16)(13,-15)(14,-15)(14,-14)(14,-13)(15,-13)(15,-12)(16,-12)
\psline(13,-16)(14,-16)(15,-16)(15,-15)(15,-14)(15,-14)(16,-14)(17,-14)(17,-13)(17,-13)(17,-12)(16,-12)
\psline(17,-12)(18,-12)(19,-12)
\psline[linewidth=0.01cm]{->}(22,-14)(17,-12)
\uput{0.5}[0](22,-14){$v$ is not a {\it bad} site}
\endpspicture
}
}
\noindent The configurations in $\Psi(A')$ have the
following three crucial properties:
\medskip
{\parindent .5cm
\item{i)} We claim that
$$
\max_{\omega'\in\Psi(A')}|\Psi^{-1}(\omega')|\leq2^{16\delta n}.\eqnum
$$
\eqlabel{\eqBornDeg}To prove \eqref{\eqBornDeg}, we first write for each $\widetilde\omega\in\Psi(A')$
$$
|\Psi^{-1}(\widetilde\omega)|\leq\sum_{I\subset\E(\L(n))}
\big|\{\omega\in\Omega_{\L(n)}:\,I_n(\omega)=I,\,\omega^{I}=
{\widetilde\omega}^{I}\}\big|\,.
$$
\noindent By \eqref{\eqInsideI} and \eqref{\eqOutsidemeasI},
the above sum contains only one term corresponding to
$I=I(\widetilde\omega)$. Hence
$$|\Psi^{-1}(\widetilde\omega)|\leq
\big|\{\omega\in\Omega_{\L(n)}:\,I_n(\omega)=I(\widetilde\omega),\,\omega^{I}={\widetilde\omega}^{I}\}\big|
\leq2^{|I_n(\widetilde\omega)|}\,,$$and the claim follows from \eqref{\eqbadedges}.
Finally, using the finite energy property and \eqref{\eqBornDeg} we get
$$\eqalign{
\Probn{A'}\leq&
\max_{\omega'\in\Psi(A')}\left|\Psi^{-1}(\omega')\right|\left(1\vee\frac{p}{1-p}\right)^{16\delta
n}\Probn{\Psi(A')}\cr
\leq&\exp(c_1\delta n)\Probn{\Psi(A')},
}\eqnum$$where $0<c_1<\infty$ is a constant.
\item{ii)} By \eqref{\eqOutsideI}, the map $\Psi$ does not modify the configuration inside
$\L(n/2)$, thus
$$
\Psi(A')\subset A.
$$

\item{iii)}By our cutting procedure we disconnect
$\Carre{h^{*}}$ from $\partial\L(n)$ hence
$$
\Psi(A')\,\subset\,
\big\{\,\L(3n/4)\nleftrightarrow\partial\L(n)\,\big\}\,.
$$
\par}\eqlabel{\eqfinitenergprop}
\noindent By the property $iii)$ and by
duality,
if the event
$\Psi(A')$ occurs,
there exists an outermost open
dual circuit $\Gamma$ in $\L(n)$ that surrounds $\L(3n/4)$.
Let $\Xi$ be the set of such dual circuits surrounding $\L(3n/4)$. For
every\def\geg{\widehat{\gamma}} $\geg\in\Xi$, we define $\Inter(\geg)$ as
the set of all the sites of $\L(n)$ that are surrounded by $\geg$ and  $\Ext(\geg)$,
the set of the sites of $\L(n)$ that are not surrounded by $\geg$ .
Note that\def\Open{\text{\rm Open}}
$$
\{\Gamma=\geg\}=\Open(\geg)\cap G_{\geg}\,,\eqnum
$$
\eqlabel{\eqcircuitext}where $\Open(\geg)=\{\forall\widehat{e}\in\geg:\
\widehat{\omega}(\widehat{e})=1\}$ and where  $G_{\geg}$ is a $\F_{\Ext(\geg)}$-measurable event.
\noindent By using  properties ii) and iii) and by \eqref{\eqcircuitext} we can write
$$\eqalign{
\Phi^{w,p}_{\L(n)}[\Psi(A')]&\leq\Phi^{w,p}_{\L(n)}[A\cap\bigcup_{\geg\in\Xi}\{\Gamma=\geg\}]\cr
&=\sum_{\geg\in\Xi}\Phi^{w,p}_{\L(n)}[A\cap G_{\geg}\vert \Open(\geg)]\ \Phi^{w,p}_{\L(n)}[\Open(\geg)].
}\eqnum$$
\eqlabel{\eqpsiaprim}
Since  $A$ is $\F_{\Inter{\geg}}$-measurable, $G_{\geg}$ is
$\F_{\Ext{\geg}}$-measurable, we can use the independence of the
$\sigma$-algebras $\F_{\Inter{\geg}}$ and $\F_{\Ext{\geg}}$ under
$\Phi^{w,p}_{\L(n)}[\,\cdot\,\vert\Open(\geg)]$ and the spatial Markov
property to get
$$\eqalign{
\Phi^{w,p}_{\L(n)}[A\cap
G_{\geg}\vert\Open(\geg)]=&\Phi^{w,p}_{\L(n)}[A\vert\Open(\geg)]\ \Phi^{w,p}_{\L(n)}[G_{\geg}\vert\Open(\geg)]\cr
=&\Phi^{f,p}_{\Inter(\geg)}[A]\ \Phi^{w,p}_{\L(n)}[G_{\geg}\vert\Open(\geg)].
}\eqnum
$$
\eqlabel{\eqcondDG}Also $A$ is an increasing event, so using \eqref{\eqcondDG}, we get
$$
\forall \geg\in\Xi\quad \Phi^{w,p}_{\L(n)}[A\cap
G_{\geg}\vert\Open(\geg)]\leq\Phi^{f,p}_{\L(n)}[A]\ \Phi^{w,p}_{\L(n)}[G_{\geg}\vert\Open(\geg)].\eqnum
$$
\eqlabel{\eqUnifBC}Using \eqref{\eqpsiaprim} and \eqref{\eqUnifBC} we obtain
 $$\eqalign{
\Phi^{w,p}_{\L(n)}[\Psi(A')]
\leq &\Phi^{f,p}_{\L(n)}[A]\sum_{\geg\in\Xi}\Phi^{w,p}_{\L(n)}[G_{\geg}\vert\Open(\geg)]\
\Phi^{w,p}_{\L(n)}[\Open(\geg)]\cr
=&\Phi^{f,p}_{\L(n)}[A]\ \Phi^{w,p}_{\L(n)}[\exists\geg\in\Xi:\ \Gamma=\geg]\leq\Phi^{p}_{\infty}[A].
}\eqnum$$
\eqlabel{\eqPsiAprimA}Combining \eqref{\eqPsiAprimA} with \eqref{\eqfinitenergprop} gives us
$$
\Phi^{w,p}_{\L(n)}[A']\leq\exp(c_1\delta n)\Phi^{p}_{\infty}[A].\eqnum
$$
\eqlabel{\eqBornAprim}Now we turn to the second term of \eqref{\eqfirstdecomp}, namely
$\smash{\Phi_{\Lambda(n)}^{w,p}[\Mn>|\L(n)|\delta\ ]}$.
Assuming that $n$ is bigger than $16\msub(\delta/2,p)/\delta$, we can apply lemma
\procref{\procbern} to get
$$
\Phi_{\Lambda(n)}^{w,p}[\Mn>|\L(n)|\delta\ ]\leq\exp\left[-\left(\frac{\delta
n}{6\msub(\delta/2,p)}\right)^2\right].\eqnum
$$
\eqlabel{\eqBornbord}Substituting \eqref{\eqBornAprim} and  \eqref{\eqBornbord} into
\eqref{\eqfirstdecomp} one has
$$
\Phi^{w,p}_{\L(n)}[A]\leq\exp(c_1\delta n)\ \Phi^{p}_{\infty}[A]+\exp\left[-\left(\frac{\delta
n}{6\msub(\delta/2,p)}\right)^2\right].\eqnum
$$
\eqlabel{\eqseconddecomp}It follows from the comments after proposition
\procref{\procLimiteJointe} that there exists a positive $\tau_c$
such that for all $p<p_c$ and $n>1$,
$$\eqalign{\Phi^{p}_{\infty}[A]\,
\leq\,|\partial\L(n/2)|\sup_{x\in\partial\L(n/2)}\Phi^{p}_{\infty}[0\leftrightarrow
x]\,
\leq\,2n\exp(-\tau_c(p_c-p)n/4).
}$$So that \eqref{\eqseconddecomp} becomes
$$
\Phi^{w,p}_{\L(n)}[A]\leq2n\exp(-\tau_c(p_c-p)n/4+c_1\delta n)+\exp\left[-\left(\frac{\delta
n}{6\msub(\delta/2,p)}\right)^2\right].\eqnum
$$
\eqlabel{\eqseconddecompbis}From \eqref{\eqseconddecompbis}, it is clear that the only way not to
destroy our estimates is to take $\delta$ at most of order
$(p_c-p)$. So let us choose $\delta={\tau_c\over 8c_1}(p_c-p)$.
Let $a>2$. By lemma \procref{\procSpeedmsub} we know that there exists a positive
constant  $c_2$ such that $\msub(\tau_c(p_c-p)/(16c_1),p)<c_2(p_c-p)^{-a}.$
Thus there exists a positive $c_3$ such that for all
$n>c_3(p_c-p)^{-1-a}$, \eqref{\eqseconddecompbis} becomes
$$
\Phi^{w,p}_{\L(n)}[A]\leq\exp(-(\tau_c/16)(p_c-p)n)+
\exp(-c_4(p_c-p)^{2+2a}n^2),\eqnum $$
\eqlabel{\eqseconddecompbisbis}where $c_4>0$.
Furthermore, we require that the first term is the main contribution, we
do this by imposing that $n>\tau_c(p_c-p)^{-1-2a}/(16c_4)$. We conclude the proof by choosing
$c=(c_3p_c^a)\vee(\tau_c/(16c_4))$ and $\lambda=\tau_c/16$.
\qed\enddemo
\noindent The last proposition permits us to control adequately the influence of boundary conditions near criticality.

\proclaim{Corollary  \procnum}Let $p\neq p_c$, $a>5$ and $\delta>0$. There exist
two positive constants $c=c(a,\delta)$ and $\lambda$ such that uniformly over the events $A\in
\F_{\L(n)}$ and uniformly over two measures $\Phi_1, \Phi_2$ in $\FK(\L(n(1+\delta)),p)$ we
have
$$
n>c|p-p_c|^{-a}\,\Rightarrow\,(1-e^{-\delta\lambda|p-p_c|n/2})^2\Phi_1[A]\leq\Phi_2[A]\leq (1+e^{-\delta\lambda|p-p_c|n/2})^2\Phi_1[A].
$$
\endproclaim\proclabel{\procWMII}
\demo{Proof}Consider $A\in\F_{\L(n)}$ and two partially wired boundary conditions $\pi_1$
and $\pi_2$ on the boundary $\partial\L((1+\delta)n)$. It is sufficient to
prove the statement for the measures $\Phi_1=\Phi_{\L((1+\delta)n)}^{\pi_1,p}$
and $\Phi_2=\Phi_{\L((1+\delta)n)}^{\pi_2,p}$. Let $m>(1+2\delta)n$  and define the  following
$\F_{\L(m)}^{\L((1+\delta)n)}$-measurable events, for $i=1,2$:
$$
W_i=\left\{
\omega\in\Omega_{\L(m)}:\eqalign{&\text{with wired boundary conditions
on $\L(m)$}\cr
&
\text{and the configuration $\omega$
on $\L(m)\setminus\L((1+\delta)n)$,}\cr
&\text{the boundary conditions
induced on $\L((1+\delta)n)$
are } \pi_i }
\right\}
$$
Since $\pi_1$ and $\pi_2$ are partially wired boundary conditions, it
is possible to find a large enough finite  $m$  such  that
$\Phi_{\L(m)}^{w,p}[W_i]>0, i=1,2$. We fix such an $m$ and write $\Phi_i[A]=\Phi_{\L(m)}^{w,p}[A\vert
W_i], i=1,2$. We note that  $d(\L(m)\setminus\L((1+\delta)n),\L(n))>\delta n/2$. Therefore, Proposition \procref{\procContConn} and an adaptation of the arguments of lemma 3.2 in \cite{\AlexIV} ensures the existence of a positive $c=c(a,\delta)$ such that
$$
n>c|p-p_c|^{-a}\quad\Rightarrow\quad|\Phi_{\L(m)}^{w,p}[A\vert
W_i]-\Phi_{\L(m)}^{w,p}[A]|\leq e^{-\delta\lambda|p-p_c|n/2}\Phi_{\L(m)}^{w,p}[A]\quad i=1,2.
$$Using the last inequality, we finally get
$$
\Phi_2[A]\geq(1-e^{-\delta\lambda|p-p_c|n/2})\Phi_{\L(m)}^{w,p}[A]\geq(1-e^{-\delta\lambda|p-p_c|n/2})^2\Phi_1[A],
$$and
$$
\Phi_2[A]\leq (1+e^{-\delta\lambda|p-p_c|n/2})\Phi_{\L(m)}^{w,p}[A]\leq (1-e^{-\delta\lambda|p-p_c|n/2})^2\Phi_1[A].
$$
\qed
\enddemo

\head\headnum{ Proof of the theorems}\endhead

\demo{Proof of Theorem \procref{\thmUR}} Since $U(\L(n))$ is
increasing, we have that
$$
\forall\Phi\in\FK(\Lt(n),p)\quad\Phi[U(\L(n))^c]\leq\Phi_{\Lt(n)}^{f,p}[U(\L(n))^c].
$$  By duality we get that
$$\eqalign{
\Phi[U(\L(n))^c]\leq&2\ \Phi_{\Lt(n)}^{f,p}[\exists\text{ an open dual path in }\Lh(n)\text{ of diameter }\geq n],
}
$$Let $a>5$. By Corollary
 \procref{\procWMII} and proposition \procref{\procLimiteJointe}
 there exist two positive constants
$c=c(a)$ and $\lambda_1$
such that for all $p>p_c$ and for all  $n>c(p-p_c)^{-a}$
we have
$$\eqalign{
\Phi_{\Lt(n)}^{f,p}[\exists&\text{ an open dual path in }\Lh(n)\text{ of
diameter}\geq n]\cr
&\leq\Phi_{\Lt(n)}^{f,p}[\exists\text{ an open dual path in }\Lh(n)\text{ of
diameter}\geq n]\cr
&\leq \Phi_{\infty}^{p}[\exists\text{ an open dual path in }\Lh(n)\text{ of
diameter}\geq n]\cr
&\leq 2n^4\exp(-\lambda_1(p-p_c)n)\cr
&\leq 2\exp((p_c-p)n(\lambda_1-4{\log n\over n(p-p_c)})),
}$$Note that there exists $n_0$ independent of everything such that
$$
\forall n>\max(n_0,c(p-p_c)^{-a})\qquad{\log n\over n(p-p_c)}\leq {n^{-1/2}\over p-p_c}\leq {1\over c}(p-p_c)^{3/2}.
$$Thus, the result follows by choosing $\lambda=\lambda_1/2$ and $c$ big enough.
To estimate the event $R$, notice that
$$
\Phi[R(\L(n),M)^{c}]\leq\Phi[U(\L(n))^{c}]+\Phi^{f,p}_{\Lt(n)}[\exists \text{ an open dual
path of diameter }\geq M].
$$Then, as before, we use Corollary \procref{\procWMII}
and proposition \procref{\procLimiteJointe} to get
$$\eqalign{
\Phi[R(\L(n),M)^{c}]&\leq\exp(-\lambda(p-p_c)n)+n^4\exp(-\lambda(p-p_c)M)\cr
&\leq(1+n^4)\exp(-\lambda(p-p_c)M).
}$$Finally, condition \eqref{\eqsizem} ensures that the prefactor does
not destroy our estimates and this concludes the proof.
\qed\enddemo
\noindent Now we turn to the estimation of the crossing cluster's size:
\demo{Proof of Theorem \procref{\thmW}}

\noindent To get \eqref{\eqIW}, one proceeds as in lemma \procref{\procbern}. For the second statement, one proceeds as in
lemma \procref{\procSpeedmsub} to prove that for every $a>9/8$, there
exists a positive constant $c=c(a,\delta)$ such that $\msup(\delta,
p)\leq C(p-p_c)^{-a}$. The desired result follows then from  \eqref{\eqIW}.
\qed\enddemo

\demo{Proof of Theorem \procref{\thmV}}
Let $\Phi\in\FK(\Lt(n),p)$. We renormalize  $\L(n)$ into $\und\L(n)$ by partitioning it into blocks $B(\und{x})$ of size $N\leq
n$ to get the renormalized box
$$
\und{\L}(n)=\{\und{x}\in\Z^2:\
(-N/2,N/2]^2+N\und{x}\subset(-n/2,n/2]^2\}.
$$Next, we  define the following events:
{\parindent=0.5cm
\item{-} For $\{\und{x}, \und{y}\}\in\E(\und\L(n))$, we denote by
$m(\und x,\und y)$ the middle point of the face between $B(\und x)$ and $B(\und y)$.  We also introduce the box $
D_{\und{x}, \und{y}}=m(\und x,\und y)+\L(\lfloor {N/4}\rfloor)
$ of width $\lfloor N/4\rfloor$ and centered at
$m(\und x, \und y)$. Then, we define
$$
K_{\und{x},\und{y}}=\{\exists\ \hbox{crossing in }
D_{\und{x}, \und{y}}\},\qquad
K_{\und x}=\bigcap_{\und{z}\in\und\L(n)\,:\,
|\und{x}-\und{z}|=1}K_{\und{x}, \und{z}}.
$$
\item{-} For $\und{x}\in\und\Lambda(n)$ and $M>0$, we define
$$\eqalign{
 R(\und{x})&=\{\exists! \hbox{
crossing cluster } C^*_{\und{x}}\hbox{ in }B(\und x)\}\cap\cr
 &\big\{
 \hbox{every open path } \gamma\subset B(\und x)\hbox{ with }
 \diam(\gamma)\geq M\hbox{ is included in }C^*_{\und{x}}
 \big\}.
}\eqnum$$
\par}
\noindent On $\und\L(n)$, we define the $0-1$ renormalized process
$(X(\und x), \und x\in\und\L(n))$ as the indicator of the occurrence of
the above mentioned events:
$$
\forall\und{x}\in\und\Lambda(n)\quad X(\und x)=\left\{\eqalign{
 1 &\hbox{ on } R(\und x)\cap K(\und x)\cr
0 &\hbox{ otherwise}} \right.
$$By Theorem \procref{\thmUR}, we get the following estimate on the
probability that a specific box is {\sl bad}. There exist
$\kappa,\lambda>0$ such that if
$$
n>N>4M>\frac{\log N}{\kappa(p-p_c)}\eqnum
$$
\eqlabel{\eqcondI}then
$$
\forall\und{x}\in\und{\L}(n)\qquad\Phi[X(\und x)=0]\leq\exp\left(-\lambda(p-p_c)M\right).\eqnum
$$
\eqlabel{\eqBadBlock}As $M$ will grow, we can restrict
ourselves to the case where there is no bad block at all and where the event $R(\L(n),N)$ is satisfied, namely
for all $\Phi\in\FK(\Lt(n),p)$, we write
$$\eqalign{
\Phi[V(\L(n),\delta)^c]\leq&\Phi[\exists\text{ a bad block }]
+\Phi[R(\L(n),N)^c]\cr
&+\Phi[\not\exists\text{ a bad block }\cap R(\L(n), N)\cap V(\L(n),\delta)^c].
}\eqnum$$
\eqlabel{\eqOneBadNoBad}

\noindent By \eqref{\eqBadBlock}, we get
$$
\Phi[\exists\text{ a bad block
}]\leq\frac{n^2}{N^2}\exp(-\lambda_1(p-p_c)M).\eqnum
$$
\eqlabel{\eqExistBad}For the second term of \eqref{\eqOneBadNoBad},we apply Theorem \procref{\thmUR} to get
$$
\Phi[R(\L(n),N)^c]\leq\exp(-\lambda_2(p-p_c) N),\eqnum
$$
\eqlabel{\eqRnot}

\noindent For the third term of \eqref{\eqOneBadNoBad}, we observe
that if there is no bad block then there is one single cluster in
the renormalized process that consists of all the blocks of
$\und{\L}(n)$. By the definition of the events associated to
$(X(\und{x}),\,\und x \in\und\L(n))$, this induces one crossing
cluster $\widetilde C^*$ of $\cup_{\und x\in\und \Lambda(n)}B(\und
x)$ that contains all the crossing clusters $C^*_{\und x}$, ${\und
x}\in{\und \L}(n)$. On the other hand, since  $R(\L(n),N)$ is
satisfied,  we have that $\widetilde C^*\subset C^*$, where $C^*$
is the crossing cluster of $\L(n)$, which is guaranteed to exists
thanks to the event $U(\L(n))$. Now, we define for every $\und
x\in\und\L(n)$ the random variables
$$Y(\und x)=N^{-2}|\{v\in B(\und x):\diam(C_v)\geq M\}|$$and observe that
$$
|C^*|<(1-\delta)\theta n^2
\quad\Rightarrow\quad
\sum_{\und x\in\und\L(n)}|C^*_{\und x}|<(1-\delta)\theta n^2\,.\eqnum
$$
\eqlabel{\eqClusBoites}Yet if $B(\und x)$ is a good box then every
cluster of $B(\und x)$ that is
of diameter larger than $M$ is included in $C^*_{\und x}$, thus using
\eqref{\eqOneBadNoBad}, \eqref{\eqExistBad},\eqref{\eqRnot}, \eqref{\eqClusBoites} and
by the FKG inequality we get
$$
\Phi[V(\L(n),\delta)^c]\leq2\frac{n^2}{N^2}\exp(-\lambda_3(p-p_c)M)+
\Phi\Big[\left(\frac{N}{n}\right)^{2} \sum_{ \und x\in\und\L(n)
}Y({\und x})\leq (1-\delta)\theta\vert E\Big],\eqnum
$$
\eqlabel{\eqBornVcI}where $E$ is the event that all the edges that
touch the boundary of the boxes $B(\und x)$ are {\sl closed} and
$\lambda_3=\min(\lambda_1,\lambda_2)$. Now we choose $N$ and $M$
such that the mean of the random variables $ Y(\und x)$ is big
enough: by using Corollary \procref{\procWMII} we have for $\und
x\in\und\L(n)$
$$\eqalign{
\Phi^{f,p}_{\L(N)}[
Y({\und x})
]&\geq
N^{-2}\Phi^{f,p}_{\L(N)}[|\{x\in\L(N-4M):x\leftrightarrow\partial\L(2M)+x\}|]\cr
&\geq
N^{-2}\sum_{x\in\L(N-4M)}\Phi^{f,p}_{\L(N)}[x\leftrightarrow\partial\L(2M)+x]\cr
&\geq
N^{-2}\sum_{x\in\L(N-4M)}\Phi^{f,p}_{x+\L(4M)}[x\leftrightarrow\partial\L(2M)+x]\cr
&\geq(1-e^{-(p-p_c)M/2})^2\frac{(N-4M)^2}{N^2}\Phi^{w,p}_{\L(4M)}[0\leftrightarrow\partial\L(2M)]\cr
&\geq(1-e^{-(p-p_c)M/2})^2\frac{(N-4M)^2}{N^2}\Phi^{p}_{\infty}[0\leftrightarrow\partial\L(2M)]\cr
&\geq (1-2e^{-(p-p_c)M/2})(1-\frac{8M}{N})\theta\,.
}$$By Onsager's formula, we have
$$\theta\,=\,(p-p_c)^{1/8}+o((p-p_c)^{1/8})\,,\qquad p\downarrow p_c \,.$$Thus if we choose
$$
M=\frac{\delta}{32}N\quad
\text{ and }\quad
M(p-p_c)\geq c,\eqnum
$$
\eqlabel{\eqcondII}where $c>0$ is a large enough constant we get
$$
\forall\und x\in\und\L(n)\qquad
\Phi^{f,p}_{\L(N)}[
Y({\und x})
]\geq\theta(1-\frac{\delta}{2})\,.\eqnum
$$
\eqlabel{\eqMoyYgrand}The random variables
$Y({\und x})$, ${\und x\in\und\L(n)}$, take their values in $[0,1]$
and they
are independent under
$\Phi[\cdot\ \vert E]$, thus we can use lemma \procref{\BornHoeff} with
\eqref{\eqMoyYgrand} to bound \eqref{\eqBornVcI} by
$$
\Phi[V(\L(n),\delta)^c]\leq2\frac{n^2}{N^2}\exp(-\lambda(p-p_c)M)+\exp\left(-\frac{\delta^2\theta^2(p)}{4}\frac{n^2}{N^2}\right).\eqnum
$$Let $a>3$ and $0<\alpha<1$. If $n^{\alpha}>(p-p_c)^{-a}$ and letting
$N=n^{\alpha}$, one gets
$$
\Phi[V(\L(n),\delta)^c]\leq2\exp(-\frac{\lambda}{32}\delta(p-p_c)n^{\alpha}+2(1-\alpha)\log
n)+\exp\left(-\frac{\delta^2\theta^2(p)}{4}n^{2-2\alpha}\right)\,.$$
\eqlabel{\eqcondIII}Also, under the above regime we have that
$(p-p_c)n^\alpha/\log n\rightarrow\infty$. Thus, by choosing
$n,N,M$ such that \eqref{\eqcondI} and \eqref{\eqcondII} are
satisfied, we obtain the desired
result.
\qed\enddemo

\head Appendix A\endhead
\noindent As promised, we give a proof of Proposition \procref{\procregpfister}
\demo{Proof of Proposition \procref{\procregpfister}}
From the Ising-FK coupling it follows that
$$
{1\over|\L(n)|}\Phi^{w,p}_{\L(n)}[M_{\L(n)}]={1\over |\L(n)|}\sum_{x\in\L(n)}\mu^{+,\beta}_{\L(n)}[\sigma(x)],
$$
where $\mu^{+,\beta}_{\L(n)}$ is the plus boundary condition Ising measure on $\{-1,+1\}^{\L(n)}$ taken at inverse temperature $\beta$.
The proof we present here is an adaptation of arguments included in \cite{\BricLeboPfister}.
An alternative way to derive the result is to use the ideas of \cite{\CCS}.
Let $n,k,l$ be three integers
larger than one. For $h>0$, we note $\mu^{+,\beta,h}_{n+k+l}$ the Ising measure on the box $\L(n+k+l)$ with boundary conditions $+$, at inverse temperature $\beta$ and where every spin in $\L(n+k+l)\setminus\L(n+k)$ is submitted to a positive field $h/\beta$.
Let $x\in\L(n)$.
The measure $\mu^{+,\beta,h}_{n+k+l}$ has the property that
$$
\lim_{h\uparrow\infty}\mu^{+,\beta,h}_{n+k+l}[\sigma(x)]=\mu^{+,\beta}_{\L(n+k)}[\sigma(x)].
$$
It is thus a sort of interpolation between the measures $\mu^{+,\beta}_{\L(n+k+l)}$ and $\mu^{+,\beta,h}_{\L(n+k)}$. Furthermore, it is easy to check that
$$
{\partial\mu^{+,\beta,h}_{n+k+l}[\sigma(x)]\over\partial h}=\sum_{y\in\L(n+k+l)\setminus\L(n+k)}\mu^{+,\beta,h}_{n+k+l}[\sig{x}\sig{y}]-\mu^{+,\beta,h}_{n+k+l}[\sig{x}]\mu^{+,\beta,h}_{n+k+l}[\sig{y}].
$$
Therefore, we have
$$\eqalign{
0&\leq\mu^{+,\beta}_{\Lambda(n+k)}[\sig{x}]-\mu^{+,\beta}_{\Lambda(n+k+l)}
[\sig{x}]=\cr
&\sum_{y\in\L(n+k+l)\setminus\L(n+k)}\int_{0}^{\infty}\mu^{+,\beta,h}_{n+k+l}[\sig{x}\sig{y}]-\mu^{+,\beta,h}_{n+k+l}[\sig{x}]\mu^{+,\beta,h}_{n+k+l}[\sig{y}]\,dh,\cr
}$$
Next, applying the Ising specific G.H.S inequality \cite{\Ellis}, we get that
$$\eqalign{
\mu^{+,\beta,h}_{n+k+l}[\sig{x}\sig{y}]-\mu^{+,\beta,h}_{n+k+l}[\sig{x}]&\mu^{+,\beta,h}_{n+k+l}[\sig{y}]\cr
&\leq\mu^{+,\beta}_\infty[\sig{x}\sig{y}]-\mu^{+,\beta}_\infty[\sig{x}]\mu^{+,\beta}_\infty[\sig{y}].
}$$
\noindent Note that the right hand side depends only on the infinite
volume measure.
On the other hand, by using Griffith's inequalities \cite{\Ellis}, we may estimate
$$
\mu^{+,\beta,h}_{n+k+l}[\sig{x}\sig{y}]-\mu^{+,\beta,h}_{n+k+l}[\sig{x}]\mu^{+,\beta,h}_{n+k+l}[\sig{y}]\leq\exp(-\lambda_1h),
$$
\noindent uniformly in $n+k+l$, $x,y$ and in $\beta$, where $\lambda_1$ is a positive constant.

\noindent Combining the two last inequalities with the magnetic field
representation of the boundary conditions, we finally obtain
$$\eqalign{
&\mu^{+,\beta}_{\L(n+k)}[\sig{x}]-\mu^{+,\beta}_{\L(n+k+l)}[\sig{x}]\leq\cr
&\int_{0}^{\infty}dh\sum_{y\in\L(n+k+l)\setminus\L(n+k)}\hskip -0,5cm\left\{(\mu^{+,\beta}_\infty[\sig{x}\sig{y}]-\mu^{+,\beta}_\infty[\sig{x}]\,\mu^{+,\beta}_\infty[\sig{y}])\wedge\exp(-\lambda_1h)\right\}.
}$$
\noindent First, let us consider the case where $0<\beta<\beta_c$. In this
situation, the explicit computation
(see \cite{\Mes}) yields
that the correlation is bounded above as follows:
$$
\mu^{+,\beta}_\infty[\sig{x}\sig{y}]-\mu^{+,\beta}_\infty[\sig{x}]\,\mu^{+,\beta}_\infty[\sig{y}]\leq\exp(-\lambda_2(\beta-\beta_c)|x-y|).\eqnum
$$
\eqlabel{\eqcalborne}Thus
$$\eqalign{
&\mu^{+,\beta}_{\L(n+k)}[\sig{x}]-\mu^{+,\beta}_{\L(n+k+l)}[\sig{x}]\cr
\leq&\int_0^\infty dh\hskip -0.8cm \sum_{y\in\L(n+k+l)\setminus\L(n+k)}
\hskip -0.5cm\exp(-(\lambda_1h\vee\lambda_2(\beta_c-\beta)|x-y|))\cr
\leq&\int_0^\infty dh\exp(-{\lambda_1\over 2}h)\hskip -0.5cm\sum_{y\in\L(n+k+l)\setminus\L(n+k)}
\hskip -0.5cm\exp(-{\lambda_2\over 2}(\beta_c-\beta)|x-y|)\cr
\leq&{2\over\lambda_1}\sum_{y\in\L(n+k+l)\setminus\L(n+k)}
\hskip -0.5cm\exp(-{\lambda_2\over 2}(\beta_c-\beta)|x-y|)\cr
\leq&{8\over\lambda_1}\exp(-{\lambda_2\over4}(\beta_c-\beta)k)\sum_{r=0}^l(n+k+r)\exp(-{\lambda_2\over4}(\beta_c-\beta)r)\cr
\leq&c_1{8\over\lambda_1}{n+k\over\beta_c-\beta}\exp(-{\lambda_2\over4}(\beta_c-\beta)k)\sum_{r=0}^l\exp(-{\lambda_2\over8}(\beta_c-\beta)r).
}$$
The last inequality has been obtained by bounding $n+k+r$ by $(n+k)(r+1)$ and by choosing  $c_1$ in such a way that
$$
\forall r\geq 0\qquad r+1\leq{c_1\over\beta_c-\beta}\exp({\lambda_2\over8}(\beta_c-\beta)(r+1)).
$$
Sending $l$ to infinity yields
$$
\mu^{+,\beta}_{\L(n+k)}[\sigma(x)]\leq{8c_1\over\lambda_1\lambda_2}{n+k\over(\beta_c-\beta)^2}(8+\lambda_2(\beta_c-\beta))\exp(-{\lambda_2\over 4}(\beta_c-\beta)k).
$$
Thus, there exists a positive constant $c_2$ such that
$$
\mu^{+,\beta}_{\L(n+k)}[\sigma(x)]\leq c_2{n+k\over(\beta_c-\beta)^2}\exp(-{\lambda_2\over 4}(\beta_c-\beta)k).
$$

\noindent Applying this inequality to the box $\Lambda(n)$ and the sites
in $\Lambda(n-k)$, we
deduce that for all $k<n$
$$
{1\over n^2}\sum_{x\in\L(n)}\mu^{+,\beta}_{\L(n)}[\sigma(x)]\leq
2{k\over n}+c_2{n\over(\beta_c-\beta)^2}\exp(-{\lambda_2\over 4}(\beta_c-\beta)k).
$$
\noindent We fix $\xi>0$ and $a>\xi+1$. For all $0<\beta<\beta_c$ we take $n>(\beta_c-\beta)^{-\xi}$ and choose $k=(\beta_c-\beta)^\xi
n$ and obtain
$$\eqalign{
{1\over
n^2}\sum_{x\in\L(n)}&\mu^{+,\beta}_{\L(n)}[\sigma(x)]\leq
\cr
2(\beta_c-\beta)^{\xi}
+&{c_2\over{(n(\beta_c-\beta)^a)}^{2\over a}}\exp\left(\big(-{\lambda_2\over 4}(\beta_c-\beta)^{1+\xi}{n\over\log n}+1+{2\over a}\big)\log n\right).
}\eqnum$$
\eqlabel{\eqpfisxi}The last expression suggests to impose a regime on $(\beta_c-\beta)^an$. Indeed, there
exists a positive $n_0=n_0(\xi,a)$ such that $n>n_0$ implies that
$n/\log n>n^{\xi+1\over a}$, hence by imposing
$$
(\beta_c-\beta)^an
>
\Big({8\over\lambda_2}(1+{2\over a})
\Big)^{a/(\xi+1)}\vee1\vee\beta_c^an_0
$$
\noindent we obtain that
$$\eqalign{
&{1\over
n^2}\sum_{x\in\L(n)}\mu^{+,\beta}_{\L(n)}[\sigma(x)]\leq2(\beta_c-\beta)^\xi+c_2\exp\left(-{\lambda_2\over 8}(\beta_c-\beta)^{\xi+1}n\right)\cr
\leq&(\beta_c-\beta)^\xi\left[2+c_2\exp\left(\log{1\over\beta_c-\beta}\Big(\xi-{\lambda_2\over 8}{(\beta_c-\beta)^{\xi+1}\over\log{1\over\beta_c-\beta}}n\Big)\right)\right],
}$$
Finally, note that there exists a positive $\eps=\eps(\xi,a)$ such
that for all $\beta_c-\eps\leq\beta\leq\beta_c$
$$
{(\beta_c-\beta)^{\xi+1}\over\log{1\over\beta_c-\beta}}>(\beta_c-\beta)^a.
$$
This implies that
$$\eqalign{
&{1\over
n^2}\sum_{x\in\L(n)}\mu^{+,\beta}_{\L(n)}[\sigma(x)]\cr
\leq&(\beta_c-\beta)^\xi\left[2+c_2\exp\left(\log{1\over\beta_c-\beta}\Big(\xi-{\lambda_2\over 8}{(\beta_c-\beta)^{\xi+1}\over\log{1\over\beta_c-\beta}}n\Big)\right)\right]\cr
\leq&(\beta_c-\beta)^\xi\left[2+c_2(\beta_c-\beta)^{-\xi+{\lambda_2\over 8}(\beta_c-\beta)^an}\right],
}$$
Thus, if we impose that
$$
n>{16\xi\over\lambda_2}(\beta_c-\beta)^{-a},
$$
we get that

$$
{1\over n^2}\sum_{x\in\L(n)}\mu_{\L(n)}^{+,\beta}[\sigma(x)]\leq(\beta_c-\beta)^\xi[2+c_2(\beta_c-\beta)^\xi].
$$
\noindent Thus, there exists a positive $c=c(\xi,a)$ such that for all
$\beta_c-\eps\leq\beta\leq\beta_c$ and for all
$n>c(\beta_c-\beta)^{-a}$ we have that
$$
{1\over n^2}\sum_{x\in\L(n)}\mu^{+,\beta}_{\L(n)}[\sigma(x)]
\leq\rho'(\beta_c-\beta)^{\xi},
\eqnum$$
\eqlabel{\eqpolydec}

\noindent where $\rho'$ is a positive constant. When  $\beta<\beta_c-\eps$,
\eqref{\eqpolydec} also holds, provided $\rho'$ is replaced by $\rho=\rho'\vee \eps^{-\xi}$.
\par\noindent In order to treat the case where $\beta>\beta_c$,
one proceeds in the same way. In this situation \eqref{\eqcalborne} is
replaced by the following bound that can be obtained from the results
of  \cite{\CoyWu}: there exist positive constants $\lambda_3,c_3$ and
$\delta$ such that for all  $x,y\in\Z^2$ satisfying
$|x-y|(\beta-\beta_c)>1/\delta$ we have that
$$
\mu^{+,\beta}_\infty
[\sigma(x)\sigma(y)]-m^{*}(\beta)^2\leq c_3\exp(-\lambda_3(\beta-\beta_c)|x-y|)
\,.
$$
\qed
\enddemo


\Refs \widestnumber\key{16}
 \ref\no\AlexI
 \by K. S. Alexander\paper On weak mixing in lattice models
 \jour Probab. Theory Relat. Fields \vol 110 \yr 1998 \pages 441-471
 \endref
 \ref\no\AlexII
 \by K. S. Alexander \paper Mixing properties and exponential decay
 for lattice systems in finite volumes
\jour Ann. Probab.
\vol 32
\yr 2004
\pages 441-487
 \endref
\ref\no\AlexIII
\by K. S. Alexander \paper Stability of the Wulff minimum and
 fluctuations in shape for large finite clusters in two-dimensional
 percolation
\jour Probab. Theory Related Fields \vol 91 \yr 1992 \pages 507-532
\endref
\ref\no\AlexIV
\by K. S. Alexander \paper Cube-root boundary fluctuations for
 droplets in random cluster models \jour Comm. Math. Phys. \vol 224
 \yr 2001 \pages 733-781
\endref
\ref\no\ACC
\by K. S. Alexander, J. T. Chayes, L. Chayes
\paper The Wulff construction and asymptotics of the finite cluster
distribution for two-dimensional Bernoulli percolation \jour Comm. Math. Phys. \vol 131 \yr 1990  \pages 1-50
\endref
\ref\no\Bar
\by D. Barbato \paper Tesi di Laurea \jour Universita di Pisa \yr 2002
\endref
\ref\no\Bodineau
\by T. Bodineau \paper The Wulff construction in three and more
dimensions \jour Comm. Math. Phys. \vol 207\yr 1999 \pages 197-229
\endref
\ref\no\BodineauII
\by  T. Bodineau
\paper   Slab percolation for the Ising model
\jour Probab. Theory Related Fields
\vol 132
\yr 2005
\pages 83-118
\endref
\ref\no\BricLeboPfister
\by J. Bricmont,  J.L. Lebowitz, C.E. Pfister
\paper On the local structure of the phase separation line in the two-dimensional Ising system.
\jour J. Statist. Phys.
\vol 26
\yr 1981
\pages  313-332
\endref
\ref\no\CamiaNewman\by F. Camia, C. M. Newman
\paper The Full Scaling Limit of Two-Dimensional Critical Percolation
\jour preprint
\yr 2005
\endref
\ref\no\Cerf
\by R. Cerf \paper Large deviations for three-dimensional
supercritical percolation \jour Ast\'erisque \vol 267 \yr 2000
\endref
\ref\no\CerfSF
\by R. Cerf \paper The Wulff crystal in Ising and percolation models
\jour Ecole d'\'et\'e de probabilit\'es, Saint Flour  \yr 2004
\endref
\ref\no\CerfMessikh
\by R. Cerf, R. J. Messikh
\paper On the Wulf crystal of the 2d-Ising model near criticality
\yr 2006
\endref
\ref\no\CerfPiszI
\by R. Cerf, \'A. Pisztora \paper On the Wulff crystal in the Ising
model
\jour Ann. Probab. \vol 28 \yr 2000 \pages 947-1017
\endref
\ref\no\CerfPiszII
\by R. Cerf, \'A. Pisztora \paper Phase coexistence in Ising, Potts and
percolation models \jour Ann. I. H. P. \vol PR 37 \yr 2001 \pages
643-724
\endref
\ref\no\CCS
\by J. T. Chayes, L. Chayes, R. H. Schonmann
\paper Exponential decay of connectivities in the two-dimensional Ising model
\jour J. Stat. Phys.
\vol 49
\yr 1987
\pages 433-445
\endref
\ref\no\CourMes
\by O. Couronn\'e, R.J. Messikh
\paper Surface order large deviations for 2d FK-percolation and Potts
 models
\jour Stoch. Proc. Appl. \vol 113 \yr 2004 \pages 81-99
\endref
\ref\no\DKS
\by R. L. Dobrushin, R. Koteck\'y, S. B. Shlosman
\paper Wulff construction: a global shape from local interaction
\jour Amer. Math. Soc. Transl. Ser. \yr 1992
\endref
\ref\no\ES
\by R.G. Edwards and A.D. Sokal
\paper Generalization of the Fortuin-Kasteleyn-Swenden-Wang
representation and Monte Carlo algorithm
\jour Phys. Rev. D
\vol 38
\pages 2009-2012
\yr 1988
\endref
\ref\no\Ellis
\by R. S. Ellis
\paper Entropy, Large Deviations, and Statistical Mechanics
\jour Springer-Verlag New York Inc.
\yr 1985
\endref
\ref\no\ForKast
\by C.M. Fortuin and R.W. Kasteleyn
\paper On the random-cluster model I. Introduction and relation to other models
\jour Physica
\vol 57
\pages 125-145
\yr 1972
\endref
\ref\no\GrimoII
\by G. R. Grimmett
\paper Percolation and disordered systems
\jour Lectures on Probability Theory and Statistics. Lectures from the 26th Summer school on Probability Theory held in Saint Flour, August 19-September 4, 1996 (P. Bertrand, ed.)  Lecture Notes in Mathematics
\vol 1665
\yr 1997
\endref
\ref\no\Grim
\by G. R. Grimmett
\paper The random cluster model
\jour in Probability, Statistics and Optimization: A Tribute to Peter Whittle
ed. F. P. Kelly
\yr 1994
\pages 49-63
\endref
\ref\no\GrimII
\by G. R. Grimmett
\paper The stochastic random-cluster process and the uniqueness of random-cluster measures
\jour Ann. Probab.
\vol 23
\yr 1995
\pages 1461-1510
\endref
\ref\no\GrimMarst
\by G. R. Grimmett, J. M. Marstrand
\paper The supercritical phase of percolation is well behaved
\jour Prc. R. Soc. Lond. Ser. A
\vol 430
\yr 1990
\pages 439-457
\endref
\ref\no\Hoeff
\by W. Hoeffding
\paper Probability Inequalities for Sums of Bounded Random Variables \jour Journal of the American Statistical
 Association \vol 58 \yr 1963 \pages 13-30
\endref
\ref\no\IoffeI \by D. Ioffe \paper Large deviation for the 2D
 Ising model: a lower bound without cluster expansions \jour
 J. Stat. Phys. \vol 74 \yr 1993 \pages 411-432
\endref
\ref\no\IoffeII \by D. Ioffe \paper Exact large deviation bounds up to
 $T_c$ for the Ising model in two dimensions \jour Probab. Theory
 Related Fields \vol 102 \yr 1995 \pages 313-330
\endref
 \ref\no\IoffSchon \by D. Ioffe, R. Schonmann \paper
 Dobrushin-Koteck\'y-Shlosman Theorem up to the critical temperature
 \jour Comm. Math. Phys. \vol 199 \yr 1998 \pages 117-167
\endref
\ref\no\LaaMessRuiz
\by L. Laanait, A. Messager, J. Ruiz
\paper Phase coexistence and surface tensions for the Potts model
\jour Comm. Math. Phys.
\vol 105
\pages 527-545
\yr 1986
\endref
\ref\no\CoyWu
\by B.M. Mc Coy, T.T. Wu
\book The Two Dimensional Ising Model
\publ Cambridge, MA: Harvard University Press \yr 1973
\endref
\ref\no\Mes
\by R.J. Messikh
\paper The surface tension of the 2d Ising model near criticality
\jour submitted
\yr 2006
\endref
\ref\no\Ons \by L. Onsager\paper Crystal statistics. I. A two-dimensional model with an order-disorder transition \jour  Phys. Rev. \vol  65,  \yr 1944 \pages 117-149
\endref
\ref\no\Pfis \by C. E. Pfister \paper Large deviations and phase
 separation in the two-dimensional Ising model \jour Helv. Phys. Acta
 \vol 64 \yr 1991 \pages 953-1054
\endref
\ref\no\PfisVelen \by C. E. Pfister, Y. Velenik \paper Large deviations
 and continuum limit in the 2D Ising model \jour Probab. Theory
 Related Fields \vol 109 \yr 1997 \pages 435-506
\endref
 \ref\no\Pisz \by \'A. Pisztora
 \paper Surface order large deviations for Ising, Potts and percolation models
 \jour Probab. Theory Relat. Fields \vol 104 \yr 1996 \pages 427-466
 \endref
\ref\no\Smir
\by S. Smirnov
\paper Critical percolation in the plane: Conformal
 invariance. Cardy's formula, scaling limits.
\jour C. R. Acad. Sci. Paris \vol 333 \pages 239-244
\yr 2001
\endref
\ref\no\Slade
\by G. Slade \paper The Lace Expansion and its Applications
\jour Ecole d'\'et\'e de probabilit\'es, Saint Flour  \yr 2004
\endref
\ref\no\SmirWer \by S. Smirnov, W. Werner
\paper  Critical exponents for two-dimensional percolation
\jour Math. Res. Lett.
\vol 8
\yr 2001
\pages  729-744
\endref
\endRefs
\enddocument